\documentclass[10pt]{article}
\usepackage{graphicx}
\usepackage{amssymb}
\usepackage{epstopdf}
\usepackage{amsmath}
\usepackage{amsfonts}

\def\phi{{\varphi}}

\DeclareSymbolFont{AMSb}{U}{msb}{m}{n}
\DeclareMathSymbol{\N}{\mathbin}{AMSb}{"4E}
\DeclareMathSymbol{\Z}{\mathbin}{AMSb}{"5A}
\DeclareMathSymbol{\R}{\mathbin}{AMSb}{"52}
\DeclareMathSymbol{\Q}{\mathbin}{AMSb}{"51}
\DeclareMathSymbol{\I}{\mathbin}{AMSb}{"49}
\DeclareMathSymbol{\C}{\mathbin}{AMSb}{"43}

\def\be{\begin{equation}}
\def\ee{\end{equation}}
\def\ber{\begin{eqnarray}}
\def\eer{\end{eqnarray}}

\def\beq{\begin{equation}}
\def\eeq{\end{equation}}

\begin{document}

\addtolength{\textheight}{0 cm} \addtolength{\hoffset}{0 cm}
\addtolength{\textwidth}{0 cm} \addtolength{\voffset}{0 cm}

\newenvironment{acknowledgement}{\noindent\textbf{Acknowledgement.}\em}{}

\setcounter{secnumdepth}{5}
 \newtheorem{proposition}{Proposition}[section]
\newtheorem{theorem}{Theorem}[section]
\newtheorem{lemma}[theorem]{Lemma}
\newtheorem{coro}[theorem]{Corollary}
\newtheorem{remark}[theorem]{Remark}
\newtheorem{claim}[theorem]{Claim}
\newtheorem{conj}[theorem]{Conjecture}
\newtheorem{definition}[theorem]{Definition}
\newtheorem{application}{Application}

\newtheorem{corollary}[theorem]{Corollary}

\title{ A characterization for   solutions of the Monge-Kantorovich mass transport problem}
\author{
Abbas Moameni
\thanks{Supported by a grant from the Natural Sciences and Engineering Research Council of Canada.}
\hspace{2mm}\\
{\it\small School of Mathematics and Statistics}\\
{\it\small Carleton University}\\
{\it\small Ottawa, ON, Canada K1S 5B6}\\
{\it\small  momeni@math.carleton.ca}\\
%
}

\maketitle

\begin{abstract}
 A measure theoretical approach is presented  to study the solutions of the Monge-Kantorovich optimal mass transport problems. This  approach together with  Kantorovich duality
provide an effective tool to  answer a long standing question about  the support of optimal plans
for  the  mass transport problem
involving  general cost functions.  We also establish
a criterion for the uniqueness.
\end{abstract}

\section{Introduction}
Let $(X, \mu)$ and $(Y, \nu)$ be two Polish probability spaces  and let $c: X \times Y  \to \R$ be a continuous function.
The Monge optimal transport problem  is to find a measurable map $T : X \to Y $ with
\[T_\# \mu=\nu \qquad \big(i.e. \, \, \nu(B)=\mu\big(T^{-1}(B)\big) \text{ for all measurable } B \subset Y \big )\]
in such a way that  $T$ minimizes the transportation cost i.e.,
\[\int_X c(x,Tx) \, d\mu=\inf_{S_\#\mu=\nu} \int_X c(x,Sx) \, d\mu. \qquad \qquad (M)\]
When a transport map $ T$  minimizes  the cost we call it an
optimal transport map.
 A relaxed version of the Monge problem  was formulated by  Kantorovich \cite{K} as   a linear optimization problem  on a convex domain. In fact,
let $\Pi(\mu,\nu)$ be the set of Borel probability measures on
$X \times Y$   which have $\mu$ and $\nu $ as marginal.
The transport cost associated to a transport plan
 $\pi \in \Pi(\mu,\nu)$ is given by
\[I_c(\pi)=\int_{X \times Y} c(x,y) \, d \pi.\]
Kantorovich's problem is to minimize
\[  \inf \big \{ I_c(\pi); \pi \in \Pi(\mu,\nu) \big \}. \qquad \qquad (MK)\]
 When a transport plan minimizes the cost, it will be called an optimal plan.
In contrary to the Monge problem, the Kantorovich problem always admits solutions as soon as the
cost function is a non-negative lower semi continuous  function (see \cite{V} for a proof). We refer to \cite{A-K-P,Br,B-D,Ca,CD,C-F-M,E-G,T-W} for existence and uniqueness
of solution to the Monge-Kantorovich problem  when the cost function is of distance form.  By now, the existence and uniqueness  is known in a wide
class of settings. Namely,  for  the general cost functions on the Euclidean space and  Manifolds, for non-decreasing strictly convex functions of the distance in
 Alexandrov spaces and
for squared distance on the Heisenberg group (see for instance \cite{A-R, Be, G, Gi,F-R, L, RR, V},  the bibliography
is not exhaustive).\\
A general criterion for existence and uniqueness of optimal transport maps known as the twist condition
  dictates  the map
\[  y \to \frac{\partial c(x,y)}{\partial x},  \]
to be  injective for fixed $x \in X.$  Under the twist condition and some regularity on the first marginal $\mu,$ the  optimal plan $\gamma$    which solves the
Monge-Kantorovich  problem $(MK)$ is supported on the graph of an optimal transport map $T,$ i.e., $\gamma=(\text{Id} \times T)_\#\mu.$ \\
 Beyond the twist condition, there is not much known about the support of an optimal plan except its numbered limb system structure when
 is extremal in the convex set $\Pi(\mu,\nu)$  and also
its  local rectifiability when the cost function is non-degenerate \cite{C-M-N, M-T-W}.\\
Our aim  is to use a measure theoretical approach   for which together with  the Kantorovich duality provide a practical   tool to study the
optimal mass transport problem. In this work,  we apply this method to characterize the support of
 optimal plans for  cost functions well beyond the twist  structure.
  The following definition is a straightforward generalization of the twist-condition.

\begin{definition} Let  $c: X \times Y \to \R$  be a function such that  $x \to c(x,y)$ is differentiable for all $y \in Y.$
\begin{itemize}
\item \large{{\bf Generalized-twist condition:}}  We say that $c$  satisfies the generalized-twist condition if
for any $x_0 \in X$ and  $y_0 \in Y$
  the set \[\Big\{y;\, \frac{\partial c(x_0,y)}{\partial x}=\frac{\partial c(x_0,y_0)}{\partial x}\Big\},\]
is a finite subset of $Y.$
\item \large{{\bf $\mathbf{m}$-twist condition:}} Let $m \in \mathbb{N}.$  We say that $c$  satisfies the $m$-twist condition if
for any $x_0 \in X$ and  $y_0 \in Y$
 the cardinality of the set \[\Big\{y;\, \frac{\partial c(x_0,y)}{\partial x}=\frac{\partial c(x_0,y_0)}{\partial x}\Big\},\]
is at most $m.$
We also say that $c$  satisfies the $m$-twist condition locally  if
for any $x_0 \in X$ and  $y_0 \in Y$ there exists a neighborhood $U$ of $y_0$ such that
 the cardinality of the set \[\Big\{y \in U;\, \frac{\partial c(x_0,y)}{\partial x}=\frac{\partial c(x_0,y_0)}{\partial x}\Big\},\]
is at most $m.$\\
\end{itemize}
\end{definition}
Note that the $m$-twist condition implies the generalized-twist condition, however, the converse is not true in general. We shall study the
support of optimal plans for this new class of cost functions.
We start with the following definition.
\begin{definition}\label{union}
 Say that a measure $\gamma \in \Pi(\mu, \nu)$ is supported on the graphs of  measurable maps  $\{T_i\}_{i=1}^k$  from $X$ to $Y$,
if there exists a sequence of measurable non-negative real functions
$\{\alpha_i\}_{i=1}^k$ from $X$ to $\R$ with $\sum_{i=1}^k \alpha_i(x)=1$ such that for each measurable set $S \subset X \times Y,$
\[\gamma(S)=\sum_{i=1}^k\int_X \alpha_i(x) \chi_S(x,T_ix) \, d\mu,\]
where $\chi_S$ is the indicator function of the set $S.$ In this case we write $\gamma=\sum_{i=1}^k \alpha_i(Id \times T_i)_\# \mu.$
\end{definition}

Here we state our main result.
\begin{theorem} \label{main} Let $X$ be a complete separable  Riemannian manifold  and  $Y$  be a Polish space equipped with  Borel probability measures $\mu$ on $X$
 and $\nu$ on $Y.$  Let   $c : X \times Y \to \R$ be a  bounded continuous cost function  and assume that:
\begin{enumerate}
     \item  the cost function $c$   satisfies  the $m$-twist condition;
  \item  $\mu$ is non-atomic and   any $c$-concave function is differentiable $\mu$-almost surely on its domain.
\end{enumerate}
  Then for each optimal plan
  $\gamma $ of  $(MK),$
 there exist $k\in \{1,...,m\},$ a sequence $\{\alpha_i\}_{i=1}^k$ of  non-negative functions from $X$ to
$[0, 1]$,
 and  Borel measurable maps $G_1,...,G_k$ from $ X$ to $ Y$ such that
\begin{eqnarray}
\gamma =\sum_{i=1}^k \alpha_i (\text{Id} \times  G_i)_\# \mu,
\end{eqnarray}
where   $\sum_{i=1}^k \alpha_i(x) =1$ for $\mu$-almost every $x \in X.$
\end{theorem}

In section (4) we shall  provide a criterion for the uniqueness of measures in $\Pi(\mu,\nu)$ that are supported on the union of the  graphs of a finite number of measurable maps.
We also have the following result for costs with the generalized twist property.

\begin{theorem}\label{mm}
 Under the assumption of Theorem \ref{main}, if one replaces the $m$-twist condition by the generalized-twist condition then  for each optimal plan
  $\gamma $ of  $(MK),$
 there exist a sequence (possibly infinite) of  non-negative real functions
$\{\alpha_i\}_{i=1}^k $ on $X$
 and a sequence of   Borel measurable maps $\{G_i\}_{i=1}^k$ from $X$ to $Y$  such that
\begin{eqnarray}
\gamma =\sum_{i=1}^k \alpha_i (\text{Id} \times  G_i)_\# \mu,
\end{eqnarray}
where  $\sum_{i=1}^k \alpha_i =1$ for $\mu$-almost every $x \in X.$
\end{theorem}

When $m=1,$ the $1$-twist condition is the well-known twist condition for which not only the minimizer of the Kantorovich problem is concentrated on
the graph of a single map but it is also unique provided the first marginal does not charge small sets.  For examples of a $2$-twist condition let us consider the  function $c: [0,2 \pi]\times[0,2\pi) \to \R$
defined by $c(x,y)=1-\cos(x-y).$ It obviously satisfies the $2$-twist condition. Assume first that $\mu=\nu$, in this case  the unique solution
 would have support on the graph $y = x$. However, the model becomes much more interesting when the densities associated
with $\mu$ and $\nu$ are different. We refer to \cite{C-M-N} and \cite{G-M2} where it is  proved that the optimal map  associated to the cost function $c(x,y)=1-\cos(x-y)$
 is unique and concentrated on the union of the graphs of two maps.\\

The most interesting examples of costs satisfying the generalized-twist condition are non-degenerate costs on smooth $n$-dimensional manifolds $X$ and $Y$.   Denote by $D^2_{xy} c(x_0,y_0)$
the $n\times n$ matrix of mixed second order partial derivatives of the function $c$ at the point $(x_0,y_0).$
A cost $c \in C^2(X \times Y)$ is non-degenerate provided
$D^2_{xy} c(x_0,y_0)$ is non-singular, that is $\text{det} \Big( D^2_{xy} c(x_0,y_0)\Big)\not=0$ for all $(x_0,y_0) \in X \times Y.$   Non-degeneracy is one of the main hypothesis
in the smoothness proof for optimal maps when the cost function satisfies the twist condition \cite{M-T-W}. It is also shown in \cite{M-P-W} that for
non-degenerate costs -not necessary twisted- the support of each optimal plan concentrates on some $n$-dimensional Lipschitz submanifold,  however,  their proof says little about
the submanifold itself. Note  that the
non-degeneracy condition will imply that the map $y \in  Y \to \frac{\partial  c(x,y)} {\partial x}$
is locally injective but not necessarily globally.  Indeed,  the  non-degeneracy property  implies that the cost function $c$ satisfies the  the $1$-twist condition locally.\\
We shall show that local $m$-twistedness implies  the generalized-twist condition and therefore
 one obtains a full characterization of the support of optimal plans for such cost functions due to Theorem \ref{mm}.
\begin{proposition}\label{non} Let $X$ and $Y$ be two smooth $n$-dimensional manifolds.  Assume that $c$ is continuously differentiable with respect to the first variable,
 and that it satisfies the $m$-twist condition locally for
 some $m \in \mathbb{N}.$
If for each $x \in X$ and $\alpha \in \mathbb{R}^+$ the set  $\{y \in Y;\,|\partial c(x,y)/\partial x |= \alpha \}$ is compact 
  then $c$ satisfies the generalized-twist condition.
\end{proposition}

To conclude the introduction, we shall emphasize that Theorem \ref{main} is in fact an effortless application of the methodology  presented in this work.

The manuscript is organized as follows: in the next section, we shall discuss the key ingredients  for our methodology in this work. In the third section we proceed with the proofs  of the main results,
while the final section is reserved to address the  uniqueness issue for $m$-twisted cost functions.
\section{Measurable weak sections and extremality}
 Let $(X,  \mathcal{B}, \mu)$ be a finite, not necessarily complete measure space, and $(Y, \Sigma)$  a measurable space.
 The completion of $\mathcal{B}$
with respect to $\mu$ is denoted by $\mathcal{B}_\mu,$ when necessary, we identify $\mu$ with its
completion on $\mathcal{B}_\mu.$ For  $A, B \in \mathcal{B},$ we write $A \subset_\mu B$ provided $\mu(A \setminus B)=0.$
 Similarly we define $A=_\mu B$ if and only if  $A \subset_\mu B$ and $B \subset_\mu A.$
  A function  $T: X \to Y$ is said to be $(\mathcal{B}, \Sigma)$-measurable if and only
if $T^{-1}(A) \in \mathcal{B}$ for all $A \in \Sigma.$
 The push forward of the measure $\mu$ by the map $T$ is denoted by
$T_\# \mu,$ i.e.
\[T_\# \mu (A)=\mu(T^{-1}(A)), \qquad \forall A \in \Sigma.\]
 By the change of variable formula it amounts to saying that $\int_{Y} f(y) \, d (T_\# \mu )=\int_X f \circ T (x)\, d\mu,$  for all bounded measurable functions $f: Y \to \R.$\\
 Two measurable functions $T, S : X \to Y$ are weakly equivalent, denoted by $S =_\mu T$,
iff $T^{-1}(A)=_\mu S^{-1}(A)$  for all $A \in \Sigma.$ We also have the following definition.
\begin{definition} Let $T: X \to Y$  be $(\mathcal{B}, \Sigma)$-measurable and $\nu$ a positive measure on $\Sigma.$\\
i. We call a map  $F: Y \to X$  a $(\Sigma_\nu,\mathcal{B})$-measurable  section of $T$ if $F$ is $(\Sigma_\nu,\mathcal{B})$- measurable
 and  $T \circ F = \text{Id}_Y.$\\
ii. We call a map  $F: Y \to X$  a $(\Sigma_\nu,\mathcal{B})$-measurable weak section of $T$ if $F$ is $(\Sigma_\nu,\mathcal{B})$- measurable
 and  $T \circ F =_\nu \text{Id}_Y.$
\end{definition}
Recall that a Polish space is a separable completely metrizable topological space. A Suslin space is the image of a Polish space under a continuous mapping.
Obviously every Polish space is a Souslin space. The following theorem ensures the existence of $(\Sigma_\nu,\mathcal{B})$-measurable  sections (\cite{Bo}, Theorem 9.1.3). This is indeed a consequence of  von Neumann's selection theorem.

\begin{theorem}. Let $X$ and $Y$ be Souslin spaces and let $T : X \to Y$ be a
Borel mapping such that $T(X) = Y$. Then, there exists a mapping $F : Y \to X$
such that $T \circ F(y)= y$ for all $y \in Y$ and  $F$ is measurable with respect to every
Borel measure on $Y$.
\end{theorem}
If $X$ is a topological space we denote by $\mathcal{B}(X)$ the set of Borel subsets in
$X.$ The space of Borel probability measures on a topological space $X$ is denoted by $\mathcal{P}(X)$.
 The following result shows that every  $(\Sigma_\nu,\mathcal{B}(X))$-measurable map has  a $(\Sigma,\mathcal{B}(X))$-measurable representation  (\cite{Bo}, Corollary 6.7.6).
\begin{proposition}\label{rep}
Let $\nu$ be a finite measure on a measurable space
$(Y,  \Sigma)$, let $X$ be a Souslin space, and let $F : Y\to X$ be a $(\Sigma_\nu, \mathcal{B}(X))$-measurable
mapping. Then, there exists a mapping
$G: Y \to X $ such that $G = F$  $\nu$-a.e. and $G^{-1}(B) \in \Sigma$ for all $ B \in \mathcal B(X).$
 \end{proposition}
  For a measurable map $T: (X, \mathcal{B}(X) ) \to (Y,\Sigma, \nu)$ denote by $\mathcal{M}( T, \nu)$  the set of all  measures $\lambda$ on $\mathcal{B}(X) $ so that $T$ pushes $\lambda$ forward to $\nu,$ i.e.
\[\mathcal{M}( T, \nu)=\{\lambda; \, T_\# \lambda=\nu\}.\]
Evidently $\mathcal{M}( T, \nu)$ is a convex set. A measure $\lambda$ is an extreme point of $\mathcal{M}( T, \nu)$ if  the identity $ \lambda= \theta \lambda_1+(1-\theta )\lambda_2$
with $\theta \in (0,1) $ and $\lambda_1, \lambda_2 \in
\mathcal{M}( T, \nu)$ imply that  $ \lambda_1=\lambda_2$.
 The set of extreme points of  $\mathcal{M}( T, \nu)$ is denoted by $ext\,\mathcal{M}( T, \nu).$\\

    We recall the following result from \cite{Gr} in
 which a characterization of the set $ext\,\mathcal{M}( T, \nu)$ is given (see also \cite{Ed} for the case where $T$ is continuous).

\begin{theorem}\label{Graf} Let $(Y,\Sigma, \nu)$ be a probability space, $(X, \mathcal{B}(X))$ be a Hausdorff space with a Radon probability
measure $\lambda$,
and let $T : X\to Y$ be an $(\mathcal{B}(X), \Sigma)$-measurable mapping.
  The following conditions are equivalent:\\
(i) $\lambda$ is an extreme point of $M(T, \nu)$;\\
(ii) there exists a $(\Sigma_\nu,\mathcal{B}(X))$-measurable weak section $F : Y \to X$ of the mapping
$T$ such that $\lambda = F_\# \nu$.\\

If $T$ is surjective and $\Sigma$ is countably separated, then conditions (i) and (ii)
are also equivalent to the following condition:\\
(iii) there exists a $(\Sigma_\nu,\mathcal{B}(X))$-measurable section $F : Y \to X$ of the mapping
$T$ with $\lambda = F_\# \nu$.\\
Finally, if in addition, $\Sigma$ is countably generated and for some $\sigma$-algebra $S$
with $\Sigma \subset S \subset \Sigma_\nu$, there exists an $(S,\mathcal{B}(X))$-measurable section of the mapping $T$,
then the indicated conditions are equivalent to the following condition:\\
(iv) there exists an $(S,\mathcal{B}(X))$-measurable section $F$ of the mapping $T$ such
that $\lambda = F_\# \nu$.\\
\end{theorem}
 The most interesting for applications is the case
where $X$ and $Y$ are Souslin spaces with their Borel $\sigma$-algebras and $T : X \to Y$
is a surjective Borel mapping. Then the conditions formulated before assertion
(iv) are fulfilled if we take for $S$ the $\sigma$-algebra generated by all Souslin
sets. Thus, in this situation, the extreme points of the set $M(T, \nu)$ are exactly the
measures of the form $F_\# \nu$, where $F : Y \to X$  is measurable with respect
to $(S,\mathcal{B}(X))$ and $T\circ F(y)
 = y$ for all $y \in  Y$.\\

 We shall now make use of the Choquet theory in the setting of noncompact sets of measures to represent each  $ \lambda \in M(T, \nu)$ as a
 Choquet type integral over $ ext \, M(T, \nu).$ Let us first recall some notations from von Weizs\"acker-Winkler \cite{WW}. In the measurable  space $(X,\mathcal{B}(X))$, let $H$ be a set of
 non-negative measures on $\mathcal{B}(X).$ By $\sum_H$ we denote the $\sigma$-algebra over $H$ generated by the functions $\varrho \to \varrho(B),$  $B \in  \mathcal{B}(X).$
If $H$ is a convex set of measures we denote by  $ext\, H$  the set of extreme points of $H.$ The set of tight positive measures on $\mathcal{B}(X)$ is denoted by $\mathcal{M}^+(X).$
For a family $\mathcal{F}$ of real valued functions on $X$ we define \[\mathcal{M}_{\mathcal{F}}^+(X)=\{\varrho \in\mathcal{M}^+(X); \, \mathcal{F} \subset \mathcal{L}^1(\varrho)\},\] and
$\sigma \mathcal{M}_{\mathcal{F}}^+(X)$ is the  topology on $\mathcal{M}_{\mathcal{F}}^+(X)$ of the functions $\varrho\mapsto \int f \, d \varrho,$ $f \in\mathcal{F}.$
The weakest topology on
$\mathcal{M}_{\mathcal{F}}^+(X)$ that makes the functions $\varrho\mapsto \int f \, d \varrho$ lower semi-continuous for all lower semi-continuous bounded functions $f$ on $X$ is denoted by
$v \mathcal{M}_{\mathcal{F}}^+(X).$
Denote by $v\sigma \mathcal{M}_{\mathcal{F}}^+(X)$ the topology generated by $\sigma \mathcal{M}_{\mathcal{F}}^+(X)$ and $v \mathcal{M}_{\mathcal{F}}^+(X).$
Here is  the main result of von Weizs\"acker-Winkler \cite{WW} regarding
the Choquet theory in the setting of noncompact sets of measures.

\begin{theorem} \label{ww} Let $\mathcal{F}$ be a countable family of real Borel functions on a topological space $X.$ Let $H$ be a convex subset of $\mathcal{M}_{\mathcal{F}}^+(X)$ such that
$\sup_{\varrho \in H}\varrho(X)< \infty.$ If $H$ is closed with respect to $v\sigma \mathcal{M}_{\mathcal{F}}^+(X)$ then for every $\lambda \in H$ there is a probability measure $\xi$
on $\sum_{ext H}$ which represents $\lambda$ in the following sense
\[\lambda(B)=\int_{ext \, H} \varrho(B)\, d\xi(\varrho),\]
 for every $B \in \mathcal{B}(X).$
\end{theorem}
 We  now  use the above theorem to represent each $ \lambda \in M(T, \nu)$ as a
 Choquet type integral over $ ext \, M(T, \nu).$

\begin{theorem} \label{cho}
Let $X$ and $Y$ be complete separable metric spaces and $\nu$ a probability measure on $\mathcal{B} (Y).$
Let $T:(X, \mathcal{B}(X)) \to (Y, \mathcal{B}(Y))$ be a surjective measurable mapping and  let $\lambda \in M(T, \nu).$ Then there exists
a  probability measure $\xi$ on $\sum_{ext \,  M(T, \nu)}$
 such that for each $B \in \mathcal{B}(X)$,
\[\lambda(B)=\int_{ext \, M(T, \nu)} \varrho(B)\, d\xi(\varrho).\]
\end{theorem}
\textbf{Proof.} Note first that any finite Borel measure on a Polish space is tight (\cite{Al}, Theorem 12.7). Let $\mathcal{A}$ be a countably family in $\mathcal{B}(Y)$ which generates $\mathcal{B}(Y)$  as a $\sigma$-field. Let
 \[\mathcal{F}=\big \{ \chi_A\circ T; \, A \in \mathcal{A} \big \},\]
where $\chi_A$ is the indicator  function of $A$. Note that $\mathcal{F}$ is  a countable family of real Borel functions on  $X.$
It is clear that $M(T, \nu)$ is closed with respect to the topology $v\sigma \mathcal{M}_{\mathcal{F}}^+(X)$. Thus, it follows from  Theorem \ref{ww}  that there exists
a  probability measure $\xi$ on $\sum_{ext \,  M(T, \nu)}$
 such that for each $B \in \mathcal{B}(X),$
\[\lambda(B)=\int_{ext \, M(T, \nu)} \varrho(B)\, d\xi(\varrho).\]
 \hfill $\square$\\
We refer the interested reader to \cite{Gr} in which a more general version of the above result is considered. Indeed, in  \cite{Gr}, S. Graf  proved that the measurable sections of $T$
 can, modulo $\nu,$ be parameterized by the pre-image measures of $\nu.$  He has also shown  that this parametrization can be done in a measurable way, i.e.  if $\tilde \Sigma$ is the
$\sigma$-field of universally  measurable subsets of $Y$ then there exists an $\sum_{ ext \, M(T, \nu)}\otimes \tilde \Sigma-\mathcal{B}(X)$ measurable map $L: ext \, M(T, \nu) \times Y \to X $
 with the following properties:\\
i. For fixed $\varrho \in ext \, M(T, \nu),$ the function $L(\varrho,.)$ is an $\tilde \Sigma-\mathcal{B}(X)$ measurable section for $T.$\\
ii. For every measurable  section $F$ for $T$ there exists $\varrho  \in ext \, M(T, \nu)$ with $L(\varrho , y)=F(y)$ for $\nu$-a.e. $y \in Y.$
It then follows that  for each bounded continuous  function $g$ on $X$ and $\lambda$ as in Theorem \ref{cho}, one has
\[\int_X g(x) \, d\lambda=\int_{ext \, M(T, \nu)}\int_X g(x)\, d\varrho(x)\, d\xi(\varrho)=\int_{ext \, M(T, \nu)}\int_Y g\big(L(\varrho,y)\big)\, d\nu(y)\, d\xi(\varrho).\]

Finally, we recall the notion of measure isomorphisms and almost homeomorphisms.
\begin{definition} Assume that $X$ and $Y$ are topological spaces with  $\mu \in \mathcal{P}(X)$ and $\nu \in \mathcal{P}(Y)$.
We
say that $(X,B(X), \mu)$ is isomorphic to $(Y,B(Y ), \nu)$ if there exists a
one-to-one map $T$ of $X$ onto $Y$ such that for all $A \in B(X)$ we have
$T(A) \in B(Y)$ and $\mu(A) = \nu\big(T(A)\big),$ and for all $B \in B(Y)$ we have
$T^{-1}(B) \in B(X)$ and $\mu\big(T^{-1}(B)\big) = \nu(B)$.\\
 We shall also say that  $(X,B(X), \mu)$ and $(Y,B(Y ), \nu)$ are almost homeomorphic  if there exists a one-
to-one Borel mapping $T$ from $X$ onto $Y$ such that $\nu = \mu\circ T^{-1} , T$ is continuous
$\mu$-a.e., and $T^{-1}$ is continuous $\nu$-a.e.
\end{definition}
The following is due to  Y. Sun \cite{Sun}.
\begin{theorem}\label{iso2}
Let $\mu$ be a Borel probability measure on a Polish
space $X$. Then the following assertions are true.\\
(i) There exist a Borel set $Y \subset [0, 1]$ and a Borel probability measure $\nu$
on $Y$ such that $(X, \mu)$ and $(Y, \nu)$ are almost homeomorphic.\\
(ii) If $\mu$ has no atoms then $(X,
\mu)$ and $([0, 1], \lambda)$, where $\lambda$ is Lebesgue
measure,  are almost homeomorphic.
\end{theorem}
 \section{Properties of optimal plans.}
 In this section we shall proceed with the proofs of  the statements in the introduction. We first state  some  preliminaries required for the proofs.
Assume that  $\gamma $ is an optimal plan for   $(MK).$
  It is standard that $\gamma \in \Pi (\mu,\nu)$ is non-atomic if and only if  at least one of $\mu$ and  $\nu$ is non-atomic (see for instance \cite{R}). Since $\mu$ is non-atomic
it follows from Theorem \ref{iso2} that the Borel measurable spaces  $(X, \mathcal{B}(X), \mu)$  and $(X \times Y, \mathcal{B}(X \times Y), \gamma)$ are  isomorphic.
Thus, there exists an isomorphism $T=(T_1, T_2)$ from $(X, \mathcal{B}(X), \mu)$  onto $(X \times Y, \mathcal{B}(X \times Y), \gamma)$. It can be easily deduced that
 $T_1 : X\to X$ and $T_2: X \to Y$ are surjective maps and
\[T_1 \# \mu=\mu \quad \& \quad T_2 \# \mu=\nu.\]
Consider the convex set \[\mathcal{M}(T_1, \mu)=\big\{ \lambda \in \mathcal{P}(X); \, T_1\# \lambda=\mu\big \},\] and note that $\mu \in \mathcal{M}(T_1, \mu).$
The following definition and proposition are  essential in the sequel.
\begin{definition}\label{gen} Denote by $\mathcal{S}(T_1)$ the set of all sections of $T_1$.
Let $\mathcal{K} \subset \mathcal{S}(T_1).$ We say that a  measurable function  $F:\big (X, \mathcal{B}(X)_\mu\big ) \to \big (X,\mathcal{B}(X)\big )$ is generated by $\mathcal{K}$ if there exist a sequence $\{F_i\}_{i=1}^\infty \subset \mathcal{K}$ such that $X=\cup_{i=1}^\infty A_i$ where 
\[A_i= \{x \in X; \, \, F(x)=F_i(x)\}.\]
We also denote by $\mathcal{G}(\mathcal{K})$ the set of all functions generated by $\mathcal{K}.$  It is easily seen that $\mathcal{K} \subseteq \mathcal{G}(\mathcal{K})\subseteq \mathcal{S}(T_1).$ 
\end{definition}

\begin{proposition} \label{vint} Let $\mathcal{K}$ be a nonempty subset of $\mathcal{S}(T_1).$ 
Then there exist $ k \in \mathbb{N} \cup\{\infty\}$ and a   sequence $\{F_i\}_{i=1}^k \subset \mathcal{G}(\mathcal{K})$  such that  the following assertions hold:
\begin{itemize}
\item[i.] For each $i\in \mathbb{N}$ with $i\leq k$ we have $\mu(B_{i})>0$  where $\{B_i\}_{i=1}^k$ is defined recursively as follows
\[B_1=X \quad \& \quad B_{i+1}=\Big \{x \in B_i; \,\, F_{i+1}(x)\not\in \{F_1(x),..., F_i(x) \}\Big\}\quad \text{provided }  k>1.\]
\item[ii.]  For all $F \in \mathcal{G}(\mathcal{K})$  we  have 
\[\mu \Big ( \Big \{x \in B_{i+1}^c\setminus B^c_{i}; \, \, F(x) \not\in \{F_1(x),..., F_i(x) \}\Big\} \Big)=0   \quad \text{provided }  k>1.\]
\item[iii.]If $k\not =\infty$ then for all $F \in \mathcal{G}(\mathcal{K})$ 
\[\mu \Big (\Big  \{x \in B_k; \, \, F(x) \not\in \{F_1(x),..., F_k(x) \}\Big\}\Big)=0.\]
\end{itemize}
Moreover, if either $k\not=\infty$ or $k=\infty$ and $\mu(\cap_{i=1}^\infty B_i)=0$ then for every $F \in \mathcal{G}(\mathcal{K})$ the measure $\varrho_F=F_\#\mu$ is absolutely continuous with respect to the measure $\sum_{i=1}^k \varrho_i $ where $\varrho_i={F_i}_\#\mu.$ 
\end{proposition}
\textbf{Proof.}
Let $F_1$ be an arbitrary function in $\mathcal{G}(\mathcal{K}).$  For  each $F \in \mathcal{G}(\mathcal{K})$ define
\[B^1_{F}:=\{x \in X; \, \, F(x)\not=F_1(x)\}.\]
If $\mu(B^1_{F})=0$ for all $F \in \mathcal{G}(\mathcal{K})$ then $k=1$ and we are done. Let us assume that there exists $\bar F_1 \in \mathcal{G}(\mathcal{K})$ with $\mu(B^1_{\bar F_1})>0.$
Let  \[G_1:=\big \{B^1_F;\, \, F \in \mathcal{G}(\mathcal{K})\big \}.\]
For $A, B \in G_1$ we say that $B \preceq A$ if and only if either  $B \subset A$ and $\mu (A\setminus B)>0$ or $A=B.$ It can be easily deduced that $\preceq$ is a partial order on $G_1.$ We shall now use the Zorn's Lemma to prove that $G_1$ has a maximal element. Let $\{B_\beta\}_{\beta \in I}$  be a chain in $G_1$ and define
\[\delta=\sup\{ \mu(B_\beta); \, \, \beta \in I\}.\]
We consider two cases. In the first case  we assume that there exists $\beta_0 \in I$ with $\mu(B_{\beta_0})=\delta$ and in the second case we assume that $\mu(B_{\beta})<\delta$ for all $\beta \in I.$\\

{\it Case I.} For this case we show that $B_{\beta_0}$ is a maximal element of $\{B_\beta\}_{\beta \in I}$. Fix $\beta \in I.$  Since $\{B_\beta\}_{\beta \in I}$ is a chain it follows that either $B_\beta \preceq B_{\beta_0}$ or $B_{\beta_0} \preceq B_{\beta}.$ We show that $B_{\beta_0} \preceq B_{\beta}$ does not occur unless $B_{\beta}=B_{\beta_0}.$ By the definition of $\delta$ we have that  $\mu(B_\beta)  \leq \delta=\mu(B_{\beta_0}).$ Thus, if $B_{\beta_0} \subseteq  B_{\beta}$  then $\mu(B_{\beta}\setminus B_{\beta_0})=\mu(B_{\beta})-\mu( B_{\beta_0})\leq 0.$ Then by the definition of $\preceq$ we must have $B_{\beta_0} = B_{\beta}.$  Therefore, if $B_{\beta_0} \not = B_{\beta}$ we must have $B_{\beta} \preceq B_{\beta_0}$ as desired. \\

{\it Case II.}  In this case there is a sequence $\{\beta_n\}$  such that $\mu(B_{\beta_{n+1}})>\mu(B_{\beta_{n}}) $ and $\lim_{n \to \infty } \mu(B_{\beta_{n}})=\delta.$
  We first show that $\cup_{\beta \in I}B_{\beta}=\cup_{n=1}^{\infty}B_{\beta_{n}}.$ In fact, otherwise there exists $\alpha \in I$ such that $B_{\alpha} \not \subseteq  \cup_{n=1}^{\infty}B_{\beta_{n}}.$  Since $\{B_\beta\}_{\beta \in I}$
is a chain it then implies that $B_{\beta_{n}} \subseteq B_\alpha$ for all $n.$ Thus, $\cup_{n=1}^{\infty}B_{\beta_{n}} \subseteq B_\alpha$ from which we obtain $\mu(B_\alpha)=\delta$  which contradicts the fact that $\mu(B_{\beta})<\delta$ for all $\beta.$   This indeed proves that $\cup_{\beta \in I}B_{\beta}=\cup_{n=1}^{\infty}B_{\beta_{n}}.$\\ 
Now note that since $\mu(B_{\beta_{n+1}})>\mu(B_{\beta_{n}})$ and $\{B_\beta\}_{\beta \in I}$ is a chain we must have  $B_{\beta_{n+1}} \supset B_{\beta_{n}}.$ Note also that  for each  $i\geq 1$  since $B_{\beta_i} \in G_1,$  there exists $F_{\beta_i} \in  \mathcal{G}(\mathcal{K})$ such that $B^1_{F_{\beta_i}}=B_{\beta_i}.$ We now define $F_0: X \to X$ by 
\begin{eqnarray*}
		F_0(x)=\left\{
		\begin{array}{lll}
			F_1(x), & x \not \in \cup_{n=1}^{\infty}B_{\beta_{n}},\\
			F_{\beta_1}(x), &  x \in  B_{\beta_1},\\
			F_{\beta_{n+1}}(x), & x \in B_{\beta_{n+1}} \setminus  B_{\beta_{n}}.
		\end{array}
		\right.
	\end{eqnarray*}
Notice that $F_0 \in \mathcal{G}(\mathcal{K})$ and therefore  $B^1_{F_0} \in G_1.$ We now verify that $B_{\beta} \preceq  B^1_{F_0}$ for all $\beta \in I$ from which the maximality of $B_{F_0}$ follows. Fix $\beta \in I.$ Since $\mu(B_\beta)<\delta$ there exists $n \in \mathbb{N}$ such that $\mu(B_\beta)< \mu(B_{\beta_n})$. Since  $\{B_\beta\}_{\beta \in I}$ is a chain we must have $B_\beta \subset B_{\beta_{n}}.$  Thus we just need to verify that $B_{\beta_n} \subset B^1_{F_0}.$ This simply follows from the fact that $B^1_{F_0}=\cup_{n=1}^{\infty}B_{\beta_{n}}.$ This completes the proof for case II.\\

Therefore $G_1$ with the partial order $\preceq$ satisfies the conditions in the Zorn's Lemma and has   a maximal element $B^1_{F_2}$ for some $F_2 \in \mathcal{G}(\mathcal{K}).$ Let 
\[B_{2}:=B^1_{F_2}=\{x \in X; \,\, F_2(x)\not=F_1(x)\}.\]
  Due to the maximality of $F_2$ we have $\mu(B_2)>0.$ Otherwise, we can define
  \begin{eqnarray*}
		\tilde F_2(x)=\left\{
		\begin{array}{ll}
			\bar F_1(x), & x \in B_2^c,\\
			F_2(x), &  x \in  B_2.
					\end{array}
		\right.
	\end{eqnarray*}
  If   $\mu(B_{2})=0$ then $\mu (B^1_{\bar F}\cap B^c_2)>0$ and since $B^1_{ F_2}\subset  B^1_{\tilde F_2}$ we must have  $B^1_{ F_2} \preceq B^1_{\tilde F_2}.$ It now follows from the maximality of $B^1_{F_2}$ that $ B^1_{F_2}=B^1_{\tilde F_2}$ which is a contradiction. This proves condition $(i)$ for $B_2.$ To prove condition $(ii),$
   take an arbitrary $F \in \mathcal{G}(\mathcal{K})$ and let \[B=\{x\in B_2^c\setminus B_1^c;\, \,  F(x)\not = F_1(x)\}.\] We show that 
$\mu(B)=0.$  Indeed, if $\mu(B)>0$ we can define $\tilde F \in \mathcal{G}(\mathcal{K})$ by 
\begin{eqnarray*}
		\tilde F(x)=\left\{
		\begin{array}{ll}
			F(x), & x \in B_2^c,\\
			F_2(x), &  x \in  B_2.
					\end{array}
		\right.
	\end{eqnarray*}
	Note that $B^1_{\tilde F} \supset B^1_{F_2}$ and $\mu(B^1_{\tilde F}) > \mu (B^1_{F_2}) $ that contradicts the maximality of $B^1_{F_2}$ in $G_1.$ This proves condition $(ii).$
	Therefore $F_1$ and $F_2$ satisfy conditions $(i)$ and $(ii).$\\
	
We can now repeat this argument to find $F_3$ (if there is any). Indeed, 	
for each $F \in  \mathcal{G}(\mathcal{K})$  we now define 
	\[B^2_{F}:=\Big \{x \in B_2; \,\,  F(x)\not\in \{F_1(x),F_2(x) \}\Big\}.\]
If $\mu(B^2_F)=0$ for all $F \in  \mathcal{G}(\mathcal{K})$ then condition $(iii)$ also holds for   $k=2$  and we are done.  Let us assume that    there exists $\bar F_2 \in  \mathcal{G}(\mathcal{K})$  with $\mu(B^2_{\bar F_{2}})>0.$
We now consider the set 
\[G_2:=\{B^2_F;\, \, F \in  \mathcal{G}(\mathcal{K})\}.\]
Following the same argument for the set $G_1,$ one can show that $G_2$ with $\preceq$ has a maximal element $B^2_{F_3}$ for some $F_3 \in \mathcal{G}(\mathcal{K}).$ We show that $F_1, F_2$ and $F_3$ satisfy conditions $(i)$ and $(ii).$
Let
\[B_3:=B^2_{F_3}=\Big \{x \in B_2; \,\,  F_3(x)\not\in \{F_1(x),F_2(x) \}\Big\}.\]
We have  that $\mu(B_3)>0.$ 
Otherwise, one  can define
  \begin{eqnarray*}
		\tilde F_3(x)=\left\{
		\begin{array}{ll}
			\bar F_2(x), & x \in B_3^c,\\
			F_3(x), &  x \in  B_3.
					\end{array}
		\right.
	\end{eqnarray*}
where $\mu(B^2_{\bar F_{2}})>0.$  If   $\mu(B_3)=0$ then $\mu (B^2_{\bar F_2}\cap B^c_3)>0$ and therefore $B^2_{F_3} \preceq B^2_{\tilde F_3} .$ It now follows from the maximality of $B^2_{F_3}$ that $ B^2_{F_3}=B^2_{\tilde F_3}$ which is a contradiction. This proves condition $(i)$ for  $B_3.$ 
 Take now an arbitrary $F \in \mathcal{G}(\mathcal{K})$ and let \[B=\big\{x\in B_3^c\setminus B_2^c:\, \,  F(x)\not \in\{F_1(x), F_2(x)\}\big \}.\] We show that 
$\mu(B)=0.$  Indeed, if $\mu(B)>0$ we can define $\tilde F \in \mathcal{G}(\mathcal{K})$ by 
\begin{eqnarray*}
		\tilde F(x)=\left\{
		\begin{array}{ll}
			F(x), & x \in B_3^c,\\
			F_3(x), &  x \in  B_3.
					\end{array}
		\right.
	\end{eqnarray*}
	Note that $B^2_{\tilde F} \supset B^2_{F_3}$ and $\mu(B^2_{\tilde F}) > \mu (B^2_{F_3}) $ that contradicts the maximality of $B_{F_3}$ in $G_3.$ This proves condition $(ii).$
	Therefore $F_1, F_2$ and $F_3$ satisfy conditions $(i)$ and $(ii).$\\

	By repeating this argument one obtains a  sequence $\{F_i\}_{i=1}^k \subset \mathcal{G}(\mathcal{K})$  satisfying $(i), (ii)$ and $(iii).$ If this sequence never stops we let $k=\infty.$
	
We shall now prove the last part of the proposition. Take an arbitrary function $F \in \mathcal{G}(\mathcal{K}).$	 We need to show that $\sigma_F=F_\#\mu$ is absolutely continiuous with respect to $\sum_{i=1}^k\sigma_i$ where $\sigma_i={F_i}_\#\mu.$
Define,
\[A_i=\{x \in X; \, \, F(x)=F_i(x)\}.\]
If $k \not=\infty$ it follows from conditions $(ii)$ and $ (iii)$ that 
\[\mu \Big(\Big   \{x \in X; \, \, F(x) \not\in \{F_1(x),..., F_k(x) \}\Big\}\Big)=0.\]
Therefore, $\mu(\cup_{i=1}^k A_i)=1 $ for  $k \not=\infty.$ Also, if $k=\infty$ and $\mu(\cap_{i=1}^\infty B_i)=0$  we have that $\mu(\cup_{i=1}^{\infty} B_i^c)=1.$ Thus it follows from condition $(ii)$ that $\mu(\cup_{i=1}^\infty A_i)=1.$ Thus, under both hypothesis in the proposition (i.e., $k \not=\infty$ or, $k=\infty$ and $\mu(\cap_{i=1}^\infty B_i)=0$) we have  $\mu(\cup_{i=1}^k A_i)=1.$
 It then implies that for every $B \in \mathcal{B}(X)$ we have, 
\begin{eqnarray*}\varrho_F(B)=\mu\big (F^{-1} (B)\big)=\mu\big (F^{-1} (B)\cap (\cup_{i=1}^k A_i)\big)&\leq& \sum_{i=1}^{k}\mu \big(F^{-1} (B) \cap A_i)\\&=&
 \sum_{i=1}^{k}\mu \big(F_i^{-1} (B) \cap A_i)\\&\leq &\sum_{i=1}^{k}\mu \big(F_{i}^{-1} (B)\big)=\sum_{i=1}^{k}\varrho_{i}(B),\end{eqnarray*}
from which the absolute continuity of $\varrho_F$ with respect to $\sum_{i=1}^k \varrho_i$ follows.
	 \hfill $\square$\\

\textbf{Completion of the proof of Theorem \ref{main}.}
Since $ \mu \in  \mathcal{M}(T_1, \mu)$,  it follows from  Theorem \ref{cho}  that   there exists a  probability measure
$\xi$ on $\sum_{ext \,  M(T_1, \mu)}$
 such that for each $B \in \mathcal{B}(X)$,

\begin{equation}\label{kre} \mu(B)=\int_{ext \, M(T_1, \mu)} \varrho(B)\, d\xi(\varrho), \qquad  \big (\varrho \to \varrho(B) \text{ is measurable}\big).\end{equation}

On the other hand, by Kantorovich duality (\cite{V}, Theorem 5.10) there exists a pair of  $c$-conjugate functions $\phi \in L^1(\mu)$ and   $\psi \in L^1(\nu)$   such that
 $\phi(x)+\psi(y) \leq c(x,y)$ for all $x,y$ and
 \begin{equation*}
 \int_{X \times Y} c(x,y) \, d\gamma= \int_X \phi(x) \, d \mu+ \int_Y \psi (y) \, d \nu.
 \end{equation*}
As $T=(T_1, T_2)$ is an isomorphism between  $(X, \mathcal{B}(X), \mu)$  and  $(X \times Y, \mathcal{B}(X \times Y), \gamma),$ it follows that

\begin{equation*}
 \int_X c(T_1 x, T_2 x) \, d\mu=\int_X \phi(T_1x) \, d \mu+ \int_X \psi(T_2 x ) \, d \mu.
 \end{equation*}
 from which together with the fact that $c(x,y) \geq \phi(x)+\psi(y)$ we obtain
 \begin{equation*} \label{fc}
 c(T_1 x, T_2 x)=\phi(T_1x) + \psi(T_2 x ).  \qquad \mu-a.e.
 \end{equation*}
 Since $\phi$ is $\mu$ almost surely differentiable and ${T_1}_{\#}\mu=\mu, $ it follows that
 \begin{equation} \label{fccc}
D_1 c(T_1 x, T_2 x)=\nabla \phi(T_1x)  \qquad \mu-a.e.
 \end{equation}
where $D_1 c $ stands for the partial derivative of $c$ with respect to the first variable.
 Let $A_\gamma \in \mathcal{B}(X)$ be the  set with $\mu(A_\gamma )=1$ such that (\ref{fccc}) holds for all $x \in A_\gamma, $ i.e.
 \begin{equation} \label{fc}
D_1 c(T_1 x, T_2 x)=\nabla \phi(T_1x)  \qquad \forall x \in A_\gamma.
 \end{equation}
Since $\mu(X\setminus A_\gamma)=0,$ it follows from (\ref{kre}) that \[\int_{ext \, M(T_1, \mu)} \varrho(X\setminus A_\gamma)\, d\xi (\varrho)=\mu(X\setminus A_\gamma)=0,\]
and therefore
there exists  a $\xi$-full measure subset $K_\gamma$ of $ext \, M(T_1, \mu)$ such that $\varrho(X\setminus A_\gamma)=0$ for all $\varrho \in K_\gamma.$ Let us now define
\[\mathcal{K}:=\big \{F\in \mathcal{S}(T_1);\,\, \exists \varrho \in K_\gamma \text{ with } \mu=F_\# \varrho   \big \},\]
where $\mathcal{S}(T_1)$ is the set of all sections of $T_1.$ Let $\mathcal{G}(\mathcal{K})$
be the set of all sections of $T_1$ generated by $\mathcal{K}$ as in Definition \ref{gen}.  By Proposition \ref{vint}, 
 there exist $k \in \mathbb{N}\cup \{\infty\} $ and a  sequence $\{F_i\}_{i=1}^k \subset \mathcal{G}(\mathcal{K})$  satisfying conditions $(i), (ii)$ and $(iii)$ in that proposition. \\

{\it Claim.} We  have  that $k \leq m.$ \\
To prove the claim assume that $k>m$ and 
for each $1\leq n\leq m+1$ let  $\varrho_n={F_n}_\# \mu.$ Since  $\varrho(X\setminus A_\gamma)=0$ for all $\varrho \in K_\gamma$  we must have $\varrho_n(X\setminus A_\gamma)=0$ for all $n.$ In fact, for a fixed $n$ since $F_n \in \mathcal{G}(\mathcal{K})$ there exists a sequence $\{F_{\sigma_i}\}_{i=1}^{\infty} \subset \mathcal{K}$ such that $X=\cup_{i=1}^{\infty}A_i$
where
\[A_i=\{x \in X; \, \, F_n(x)=F_{\sigma_i}\}.\] 
Let  $\sigma_i \in K_{\gamma}$ be such that  the map $F_{\sigma_i}$
is a push-forward from  $\sigma_i$  to $\mu.$ It follows that 
\begin{eqnarray*}
\varrho_n(X\setminus A_\gamma)=\mu\big (F_n^{-1}(X\setminus A_{\gamma})\big )&=&\mu\big ((\cup_{i=1}^{\infty}A_i) \cap  F_n^{-1}(X\setminus A_{\gamma})\big )\\
&\leq &\sum_{i=1}^\infty \mu\big (A_i \cap  F_n^{-1}(X\setminus A_{\gamma})\big )\\
&=&\sum_{i=1}^\infty \mu\big (A_i \cap  F_{\sigma_i}^{-1}(X\setminus A_{\gamma})\big )\\
&\leq &\sum_{i=1}^\infty \mu\big ( F_{\sigma_i}^{-1}(X\setminus A_{\gamma})\big )=\sum_{i=1}^\infty \sigma_i(X\setminus A_{\gamma})=0.\\
\end{eqnarray*}
This proves that for each $n \in \{1,...,m+1\}$ we have that $\varrho_n(X\setminus A_\gamma)=0.$
It now follows  from (\ref{fc}) that
 \begin{equation}
D_1 c\big (T_1\circ F_n (x), T_2 \circ F_n(x)\big)=\nabla \phi \big(T_1\circ F_n(x)\big)  \qquad  \forall x \in F_n^{-1}(A_\gamma).
 \end{equation}
Thus,
\begin{equation}\label{fc3}
D_1 c\big (x, T_2 \circ  F_n (x) \big)=\nabla \phi(x)  \qquad  \,\,\forall x \in \cap_{n=1}^{m+1}F_n^{-1}(A_\gamma),  \,\,  \forall n \in\{1,2,...,m+1\}.
 \end{equation}
 Since $\varrho_n(X \setminus A_\gamma)=0$ and $\varrho_n$ is a probability measure we have  that $\varrho_n(A_\gamma)=1$ for every $n \in\{1,2,...,m+1\}.$
Therefore, $\mu\big(F_n^{-1}(A_\gamma)\big)=\varrho_n(A_\gamma)=1$. This implies that  $\mu \big(\cap_{n=1}^{m+1}F_n^{-1}(A_\gamma)\big)=1.$   This together with (\ref{fc3}) yield that
  \begin{equation}\label{ku}
D_1 c\big (x, T_2 \circ  F_n(x)\big)=\nabla \phi(x)  \qquad  \,\,\forall x \in \bar A_\gamma,
 \end{equation}
 where $\bar A_\gamma=\cap_{n=1}^{m+1}F_n^{-1}(A_\gamma).$ Note that by condition $(i)$ in Proposition \ref{vint} we have  $\mu(B_{m+1})>0.$
 Take $x \in \bar A_\gamma\cap B_{m+1}.$ It follows from  the $m$-twist condition that the cardinality of the set 
 \[L_x:=\Big \{y \in Y; \, D_1c\big (x, T_2 \circ  F_1(x)\big )=D_1c\big (x, y\big ) \Big\},\]
 is at most $m$. On the other hand it follows from (\ref{ku}) that $T_2 \circ F_n(x) \in L_x$ for all $n \in \{1,2,...,m+1\}.$ Thus, there exist $i, j \in \{1,2,...,m+1\}$ with $i <j$ such that $T_2 \circ F_i(x)=T_2 \circ F_j(x).$
 Since  $T_1 \circ  F_i=T_1 \circ F_j=Id_X$ and  the map $T=(T_1,T_2)$ is injective it follows that $F_{i}(x)=F_j(x).$ On the other hand  $x \in B_{m+1} \subseteq B_j$ from which we have $F_j(x)\not \in \{F_1(x),...,F_{j-1}(x)\}.$ This leads to a contradiction and the claim follows.\\
 
  By the latter  claim  we have that $k\leq m.$ It now  follows from Proposition \ref{vint} that  every 
 $\varrho \in K_\gamma$ is absolutely continuous with respect to the measure $\sum_{i=1}^k \varrho_i  $ where $\varrho_i={F_i}_\# \mu$ for $1 \leq i \leq k.$ 
 This together with the representation 
 \begin{equation*} \mu(B)=\int_{ext \, M(T_1, \mu)} \varrho(B)\, d\xi(\varrho)=
 \int_{K_{\gamma}} \varrho(B)\, d\xi(\varrho), \qquad  \big (\forall B \in \mathcal{B}(X)\big),\end{equation*}
 imply that $\mu$ is absolutely continuous with respect to $\sum_{i=1}^k \varrho_i.$  
   It then follows that there exists a
non-negative measurable function $\alpha: X \to \R\cup\{+\infty\}$ such that
\[\frac{d \mu}{d \big(\sum_{i=1}^k \varrho_i\big)}=\alpha.\]
  Define $\alpha_i=\alpha \circ F_i$ for $i=1,...,k.$  We show that
 $ \sum_{i=1}^k\alpha_i(x)=1$ for $\mu$-almost every $x \in X.$ In fact, for each $B \in \mathcal{B}(X)$ we have
\begin{eqnarray*}\mu(B)=\mu(T_1^{-1}(B))=\sum_{i=1}^k\int_{T_1^{-1}(B)} \alpha(x)  \, d\varrho_i =\sum_{i=1}^k\int_{F_i^{-1}\circ T_1^{-1}(B)} \alpha(F_ix)  \, d \mu
=\sum_{i=1}^k\int_{B} \alpha_i(x)
  \, d \mu,
\end{eqnarray*}
from which we obtain $\mu(B)=\sum_{i=1}^k\int_{B} \alpha_i(x)
  \, d \mu.$ Since this holds for all $B \in \mathcal{B}(X)$ we have
\[ \sum_{i=1}^k\alpha_i(x)=1, \qquad \quad \mu-a.e.\]
  It follows from Proposition \ref{rep} that each $F_i$ is $\mu$-a.e. equal to  a $(\mathcal{B}(X),\mathcal{B}(X))$-measurable
 function for which we still denote it by $F_i.$  For each $i \in \{1,...,k\},$ let  $G_i=T_2 \circ F_i.$  We now show that
 $\gamma =\sum_{i=1}^k \alpha_i (\text{Id}\times  G_i)_\# \mu$.  For each bounded continuous function $f : X \times Y \to \R$ it follows that
\begin{eqnarray*}
\int_{X \times Y} f(x,y) \, d \gamma=\int_X f(T_1x,T_2x) \, d\mu&=&\sum_{i=1}^k\int_X  \alpha(x)f(T_1x,T_2x) \, d\varrho_i\\
&=& \sum_{i=1}^k\int_X \alpha\big (F_i (x)\big ) f\big (T_1\circ F_i (x),T_2 \circ F_i (x) \big ) \, d\mu\\
&=& \sum_{i=1}^k\int_X \alpha_i(x) f\big (x,G_i (x) \big ) \, d\mu.
\end{eqnarray*}
Therefore,
\[\gamma =\sum_{i=1}^k \alpha_i (\text{Id}\times   G_i)_\# \mu.\]
\hfill $\square$\\

 Note  that one can weaken
the assumptions on the cost function in Theorem \ref{main}. Since,  this does not require  new ideas we do not elaborate.\\

\textbf{Proof of Theorem \ref{mm}}.
 The proof goes almost in the same lines as the proof of Theorem \ref{main}. Indeed,   there exists a  probability measure
$\xi$ on $\sum_{ext \,  M(T_1, \mu)}$
 such that for each $B \in \mathcal{B}(X)$,

\begin{equation}\label{okre} \mu(B)=\int_{ext \, M(T_1, \mu)} \varrho(B)\, d\xi(\varrho), \qquad  \big (\varrho \to \varrho(B) \text{ is measurable}\big).\end{equation}
 There also exists a set  $A_\gamma \in \mathcal{B}(X)$  with $\mu(A_\gamma )=1$ such that (\ref{fccc}) holds for all $x \in A_\gamma, $ i.e.
 \begin{equation} \label{ofc}
D_1 c(T_1 x, T_2 x)=\nabla \phi(T_1x)  \qquad \forall x \in A_\gamma.
 \end{equation}
Since $\mu(X\setminus A_\gamma)=0,$ it follows from (\ref{okre}) that \[\int_{ext \, M(T_1, \mu)} \varrho(X\setminus A_\gamma)\, d\xi (\varrho)=\mu(X\setminus A_\gamma)=0,\]
and therefore
there exists  a $\xi$-full measure subset $K_\gamma$ of $ext \, M(T_1, \mu)$ such that $\varrho(X\setminus A_\gamma)=0$ for all $\varrho \in K_\gamma.$ Define
\[\mathcal{K}:=\big \{F\in \mathcal{S}(T_1);\,\, \exists \varrho \in K_\gamma \text{ with } \mu=F_\# \varrho   \big \},\]
and let $\mathcal{G}(\mathcal{K})$
be the set of all sections of $T_1$ generated by $\mathcal{K}$ as in Definition \ref{gen}.  By Proposition \ref{vint}, 
 there exist a  sequence $\{F_i\}_{i=1}^k \subset \mathcal{G}(\mathcal{K})$ and a sequence $\{B_i\}$ of 
 subsets of $X$ with  $ k \in \mathbb{N} \cup\{\infty\}$ satisfying conditions $(i), (ii)$ and $(iii)$ in that proposition. If $k\not=\infty$ then by the same argument as in the proof of Theorem \ref{main} the measure $\gamma$ is supported on the graphs of $k$ maps as desired. Let us now examine the case where $k=\infty.$
{\it Claim.} If  $k=\infty $ then   $\mu(\cap_{i=1}^\infty B_i)=0.$ \\
To prove the claim assume that $ k=\infty$ and $\mu(\cap_{i=1}^\infty B_i)>0.$
 Set $\varrho_i={F_i}_\# \mu$   for each $i\in \mathbb{N}.$  Since  $\varrho(X\setminus A_\gamma)=0$ for all $\varrho \in K_\gamma$  we must have $\varrho_i(X\setminus A_\gamma)=0$ for all $i.$  It now follows  from (\ref{ofc}) that
 \begin{equation}
D_1 c\big (T_1\circ F_i (x), T_2 \circ F_i(x)\big)=\nabla \phi \big(T_1\circ F_i(x)\big)  \qquad  \forall x \in F_i^{-1}(A_\gamma).
 \end{equation}
Thus,
\begin{equation}\label{ofc3}
D_1 c\big (x, T_2 \circ  F_i (x) \big)=\nabla \phi(x)  \qquad  \,\,\forall x \in \cap_{i=1}^{\infty}F_n^{-1}(A_\gamma),  \,\,  \forall i.
 \end{equation}
 Since $\varrho_i(X \setminus A_\gamma)=0$  we have  that $\varrho_i(A_\gamma)=1.$ 
Note that $\varrho_i(A_\gamma)=\mu\big(F_i^{-1}(A_\gamma)\big)$ and therefore $\mu \big(\cap_{i=1}^{\infty}F_i^{-1}(A_\gamma)\big)=1.$   This together with (\ref{ofc3}) yield that

  \begin{equation}\label{oku}
D_1 c\big (x, T_2 \circ  F_i(x)\big)=\nabla \phi(x)  \qquad  \,\,\forall x \in \bar A_\gamma,
 \end{equation}
 where $\bar A_\gamma=\cap_{i=1}^{\infty}F_i^{-1}(A_\gamma).$ Note that we have assumed that $\mu(\cap_{i=1}^\infty B_i)>0.$ 
 Take $x \in  (\cap_{i=1}^\infty B_i) \cap \bar A_\gamma.$ It follows from  the generalized-twist condition that the following set 
 \[L_x:=\Big \{y \in Y; \, D_1c\big (x, T_2 \circ  F_1(x)\big )=D_1c\big (x, y\big ) \Big\},\]
 is a finite subset of $Y$. On the other hand it follows from (\ref{oku}) that $T_2 \circ F_i(x) \in L_x$ for all $ i \in \mathbb{N}.$   Thus, there exist $i, j$ with $i <j$ such that $T_2 \circ F_i(x)=T_2 \circ F_j(x).$
 Since  $T_1 \circ  F_i=T_1 \circ F_j=Id_X$ and  the map $T=(T_1,T_2)$ is injective it follows that $F_{i}(x)=F_j(x).$ On the other hand  $x \in \cap_{i=1}^\infty B_i \subseteq B_j$ from which we have $F_j(x)\not \in \{F_1(x),...,F_{j-1}(x)\}.$ This leads to a contradiction and the claim follows.
 The rest of the proof is now similar to the proof of Theorem \ref{main}.
  \hfill $\square$\\

We conclude this section by proving the generalized-twist property for locally m-twisted costs.\\

\textbf{Proof of Proposition \ref{non}}. Fix $x_0 \in X$ and $y_0 \in Y.$ We need to show that
the set \[L_{(x_0,y_0)}=\Big\{y \in Y; \, D_1 c(x_0,y_0)=D_1 c(x_0,y)\Big\},\]
is finite. If $L_{(x_0,y_0)}$ is not finite there exists an infinitely countable subset $\{y_n\}_{n \in \mathbb{N}} \subset L_{(x_0,y_0)}.$  By the compactness assumption in the statement of the
 proposition \ref{non}  the sequence  $\{y_n\}_{n \in \mathbb{N}} $ has an
accumulation point $\bar y \in Y$ and there exists a subsequence still denoted by $\{y_n\}_{n \in \mathbb{N}}$ such that $y_n \to \bar y.$  Since
 $D_1 c$ is continuous it follows that
$\bar y \in L_{(x_0,y_0)}.$ Since $c$ is $m$-twisted locally, there exists a neighborhood  $U$ of $\bar y$ such that the cardinality of the set
\[\big \{y \in U;  D_1c(x_0,y)=D_1c(x_0,\bar y)\big \}\]
 is at most $m.$
This is a contradiction as $\bar y$ is an accumulation point of the
sequence $\{y_n\}$ and \[D_1 c(x_0,\bar y)=D_1 c(x_0,y_0)=D_1 c(x_0, y_n), \qquad \forall n \in \mathbb{N}.\] This completes the proof. \hfill $\square$
\section{Uniqueness}
In this section we shall discuss some cases where we may have uniqueness in Theorem \ref{main}.  Note first that uniqueness under the $1$-twist  condition has been extensively
 studied and
it is known that  if $\mu$ is absolutely continuous with respect to the volume measure then uniqueness occurs. However, this can fail if $\mu$ charges small sets.
This observation makes it evident that under the $m$-twist condition on the $x$ variable and $1$-twist condition with respect to the $y$ variable ($x \to c_y(x,y)$ is injective)
the uniqueness occurs
provided both $\mu$
and $\nu$ do not charge small sets.  There is also another uniqueness criterion  known as the sub-twist property \cite{C-M-N} i.e. for each  $y_1 \not= y_2 \in Y$ the
map $x \to c(x,y_1)-c(x,y_2)$  has no critical points, save at most one global maximum and one global minimum.
Our approach is  to study the extremality of the transport plans supported on the union of the graphs of finitely many functions in the convex set $\Pi(\mu, \nu).$

The following result shows that uniqueness may occur up to the support of  optimal plans.
\begin{proposition} Suppose that  $c$  satisfies the $m$-twist condition and all the assumptions of Theorem \ref{main} are fulfilled. Let
$\bar \gamma$ be an optimal plan such that
\begin{eqnarray}
\bar \gamma =\sum_{k=1}^m \alpha_k (\text{Id}\times  G_k)_\# \mu, \qquad  \big (\alpha_i(x) \geq 0  \text{ and } \, \alpha_1(x) \alpha_2(x)...\alpha_m(x)\not=0 \text{ for }  \mu-\text{a.e. } \, x \in X\big ),
\end{eqnarray}
and for each $i\not=j$ the set $\{x; G_i(x)=G_j(x)\}$ is $\mu$-negligible.
  Then for any other optimal plan $\gamma$  we have
 \[\text{Supp} (\gamma) \subseteq \text{Supp}(\bar \gamma).\]
\end{proposition}
\textbf{Proof.}
Take $\phi$ and $\psi$ as in the proof of Theorem \ref{main}. It follows that
\[ \int c(x,y) \, d \bar  \gamma= \int [\phi(x)+\psi(y)]\, d\bar \gamma,  \]
from which we obtain
\[\sum_{k=1}^m \int_X \alpha_i c(x,G_i x) \, d \mu=\sum_{k=1}^m \int_X \alpha_i \big[ \phi(x)+\psi (G_i x)\big] \, d \mu.\]
It then follows that
\[\sum_{k=1}^m \int_X \alpha_i [c(x,G_i x)-\phi(x)-\psi (G_i x)] \, d \mu=0.\]
Since each integrand in the latter expression is non-negative  it yields that
\[c(x,G_i x)=\phi(x)+\psi (G_i x)\qquad \mu-a.e. \quad \forall i \in\{1,...,m\}.\]
Consequently we  obtain,
\begin{equation}\label{uniq}D_1 c(x,G_i x)=\nabla \phi(x)\qquad \mu-a.e. \quad \forall i \in\{1,...,m\}.\end{equation}
Note also that  for  $i\not=j$ the set $\{x \in X; \, \, G_i(x)=G_j(x)\}$ is a null set with respect to the measure $\mu.$ This together with (\ref{uniq}) and
 the $m$-twist
condition imply that  the cardinality of
 the set $\{G_1x,...,G_m x\}$ is $m$ for $\mu$-a.e. $x \in X.$\\
 Now assume that  $\gamma$ is also  an optimal plan
   of  $(MK).$   It follows from Theorem \ref{main} that
 there exist  a sequence of non-negative functions
$\{\beta_i\}_{i=1}^m $
 and,  Borel measurable maps $T_1,...,T_m: X \to Y$ such that
\begin{eqnarray*}
\gamma =\sum_{i=1}^m \beta_i (\text{Id} \times  T_i)_\# \mu
\end{eqnarray*}
By a similar argument as above one obtains
\begin{equation*}\beta_i(x)\big [D_1 c(x,T_i x)-\nabla \phi(x)\big]=0\qquad \mu-a.e. \quad \forall i \in\{1,...,m\}.\end{equation*}
For each $i$ define $\Omega_i=\{x \in X; \beta_i(x)\not=0\}.$
Since the cardinality of the set $\{G_1x,...,G_m x\}$ is $m$ for $\mu$-a.e. $x \in X$ and since  $c$ satisfies the $m$-twist condition
we have  that for each $i,$ $\{T_i x\} \subseteq \{G_1x,...,G_m x\}$ for $\mu$-a.e. $x \in \Omega_i$. This completes the proof. \hfill $\square$\\

We shall know provide a criterion for the uniqueness of measures in $\Pi(\mu,\nu)$ that are supported on the graphs of a finite number of measurable maps. For a map $f$ from a set $X$ to a set $Y$  denote by  $Dom(f)$ the domain of $f$ and by $R(f)$ the range of $f.$  Here is our main theorem addressing the uniqueness issue. 

\begin{theorem}\label{aper} Let $X$ and  $Y$  be  Polish spaces equipped with  Borel probability measures $\mu$ on $X$
 and $\nu$ on $Y,$ and    let  $\{T_i\}_{i=1}^k$ be a sequence of measurable maps from $ X$ to $Y$.    
Assume that   the following assertions hold:
\begin{enumerate}
\item For each $i \in \{2,...,k\} $ the map  $T_i$ is injective and $R(T_i)\cap R(T_j)=\emptyset$
for all $2\leq i,j \leq k$ with $i \not=j.$
\item  There exists a bounded measurable function $\theta:Y \to \R$ with the property that $\theta(T_1x)-\theta(T_ix)\geq 0$ on $X$
 and  $\theta(T_1x)-\theta(T_ix)=0$  if and only if  $T_1x =T_ix$.
\end{enumerate}
Then there exists at most one $\gamma \in \Pi(\mu, \nu)$ that is supported on the union of  the graphs of $T_1, T_2,..., T_k.$
 \end{theorem}

 As an immediate consequence of the latter Theorem  we recover the following  uniqueness result due  to Seethoff and  Shiflett  \cite{S-S}.
\begin{corollary}
 Let $X=Y=[0,1]$ and $\mu=\nu$ is the Lebesgue measure. If $T_1 \leq T_2$ and one of $T_1$ or $T_2$ is injective on $D=\{x;  T_1(x) \not=T_2(x)\}$ then
 there exists at most one $\gamma \in \Pi(\mu, \nu)$ that is supported on the graphs of $T_1$ and $ T_2.$
\end{corollary}
\textbf{Proof} Suppose $T_2$ is injective on $D.$ One can define $\theta: Y \to \R$ by  $\theta (y)=-y.$ Since $T_1 \leq T_2$ then
$\theta \big (T_1(y)\big ) -\theta \big(T_2(y) \big)$ is non-negative  and $\theta\big (T_1(y)\big ) -\theta \big(T_2(y) \big)=0$ iff $T_1(y)=T_2(y).$ The result then follows from Theorem
\ref{aper} and Remark \ref{aper2}. \hfill $\square$\\

Here is another application of Theorem (\ref{aper}) for maps with disjoint ranges.
\begin{corollary}
 Let $X$ and  $Y$  be  Polish spaces equipped with  Borel probability measures $\mu$ on $X$
 and $\nu$ on $Y.$   Let  $\{T_i\}_{i=1}^k$ be a sequence of measurable maps from $ X$ to $Y$ such that $T_i$ is injective for each $i \in \{2,...,k\}$ and
 $R(T_i)\cap R(T_j)=\emptyset$
for all $1\leq i,j \leq k$ with $i\not=j.$ If $R(T_1)$ is measurable  then there exists at most one $\gamma \in \Pi(\mu, \nu)$ that is supported on the graphs of $T_1, T_2,..., T_k.$
\end{corollary}
\textbf{Proof.}  By considering
$\theta(y)=\chi_{R(T_1)} (y),$ the indicator function of $R(T_1),$  one can easily check that $\theta$ satisfies the required properties in assumption (2) of Theorem \ref{aper}.
\hfill $\square$\\

We need some preliminaries before proving  Theorem \ref{aper}.
  For a map $f$ from a set $X$ to a set $Y$  the graph of $f$  is denoted by $Graph (f)$ and defined by
\begin{eqnarray*}
 Graph(f)=\{(x, f(x)); \, x \in Dom(f)\}.
\end{eqnarray*}
For a map $g$ from $Y$ to $X$, the antigraph of $g$ is denoted by $Antigraph(g)$  and  defined by
\begin{eqnarray*}
 Antigraph(g)=\{(g(y),y); \, y \in Dom(g)\}.
\end{eqnarray*}
For any family of maps $F_1,...,F_n$ from a set $U$ to itself we denote by
 $O_{m=1}^n F_i$ the following operation,
\begin{eqnarray}\label{oo}
O_{m=1}^n F_i=F_1 \circ F_2 \circ ...\circ F_n,
\end{eqnarray}
where the symbol $\circ$ is simply the composition operator between two functions. 
Here we recall the definition of aperiodic representations \cite{B-S}.
 \begin{definition}\label{ap} Let $X$ and $Y$ be two sets and let $f: X \to Y$ and $g:Y \to X$. Define
\begin{eqnarray*}
		T(x)=\left\{
		\begin{array}{ll}
			g \circ f(x), & x \in Dom(f)\cap f^{-1} \big ( Dom (g) \big)=D(T),\\
			x, &  x \notin D(T).\\
		\end{array}
		\right.
	\end{eqnarray*}
The maps $f, g$ are aperiodic if $x \in D(T)$ implies that $T^n(x)\not=x$ for any $n \geq1.$\\
If  $S=Graph(f) \cup Antigraph(g), $ $Graph(f) \cap Antigraph(g)=\emptyset$ and $f,g$ are aperiodic, then this is called an aperiodic decomposition of $S.$
Moreover, if $(X, \Sigma (X))$ and $(Y, \Sigma(Y))$ are measure spaces and the maps $f$  and $g$ are measurable we call the  maps  $f, g$  measure-aperiodic  if any $T$-invariant probability measure defined on $\Sigma(X)$ is supported by
 $X\setminus D(T).$
   \end{definition}
It what follows we say that
 $\gamma \in \Pi(\mu, \nu)$ is concentrated on a set $S$ if the outer measure of its complement is zero, i.e. $\gamma^*(S^c)=0.$
 We recall the following result from \cite{B-S} regarding doubly stochastic measures with aperiodic supports.
\begin{theorem}[Benes \& Stepan 1987]\label{BS}
 Let $(X, \mathcal{B}(X), \mu)$ and $(Y, \mathcal{B}(Y), \nu)$ be complete separable Borel metric spaces. Let $f: X \to Y$ and $g:Y \to X$ be aperiodic measurable maps  and
 $Graph(f) \cap Antigraph(g)=\emptyset.$ Then there exists at most one $\gamma \in \Pi(\mu, \nu)$ that is supported on $S=Graph(f) \cup Antigraph(g)$ provided $f$ and $g$
 are measure-aperiodic.
\end{theorem}
Beginning with the work of Lindenstrauss and Douglas \cite{LD1, LD2}, Hestir and Williams \cite{HW} provided an alternate proof of the latter Theorem while further refining the structure  these graphs should
take, and rewriting them in terms of  measurable limb numbering systems.   Chiappori, McCann and Nesheim \cite{C-M-N} further improved the result of  Hestir and Williams by weakening the measurability
requirement.  For our purpose in this paper the result of Benes \& Stepan seems to be more suitable. \\

\textbf{Proof of Theorem \ref{aper}.} For each $i\geq 2,$ since $T_i$ is injective  we have that $R(T_i)$ is a measurable subset of $Y.$
Define \[g:Dom(g)= \cup_{i=2}^kR(T_i)\subset Y \to X,\] by $g(y)=T_i^{-1}(y)$ for
$y \in R(T_i)$ and note that $g$ is measurable.  Define $f: Dom(f) \to Y$ by $f(x)=T_1(x)$ where \[Dom(f)=\{x \in Dom(T_1); \, \,  T_1(x) \not=T_i(x) \text{   for all  } 2\leq i\leq k \}.\]
Note that $Graph(f) \cap Antigraph(g)=\emptyset.$ In fact, if $Graph(f) \cap Antigraph(g)\not=\emptyset$ then there exists $x \in Dom(f)$ and $y \in Dom(g)$ with
$(x,f(x))=(g(y),y).$ It then follows that $y=f(x)=T_1(x)$ and $x=T^{-1}_i(y)$ for some $2\leq i\leq k.$  This is a contradiction as $T_1(x)\not=T_i(x)$ on $Dom(f).$
 Define $T: X \to X$ as in Definition \ref{ap}, i.e., 
\begin{eqnarray*}
		T(x)=\left\{
		\begin{array}{ll}
			g \circ f(x), & x \in Dom(f)\cap f^{-1} \big ( Dom (g) \big)=D(T),\\
			x, &  x \notin D(T).\\
		\end{array}
		\right.
	\end{eqnarray*}

We shall now proceed with the rest of the proof in two steps. In the first step we show that $f$ and $g$ are aperiodic and in the second step we show that $f$ and $g$ are
measure-aperiodic.  Then  the result follows from  Theorem \ref{BS}.\\

{\it Step 1:}   Assume that  there exist $x \in D(T)=Dom(f) \cap f^{-1}\big(Dom(g) \big)$
and $n \in \mathbb{N}$ such that $(g\circ f)^n(x)=x.$ It follows from the construction of $f$ and $g$ that 
\[x=(g\circ f)^n(x)=(T_{i_1}^{-1} \circ T_1 \circ T_{i_2}^{-1}\circ T_1\circ...\circ T_{i_n}^{-1} \circ T_1)(x),\]
for some $i_1,...,i_n \in \{2,...,k\}.$ As in (\ref{oo}) we will rewrite the latter expression in a short form by  \[O_{m=1}^n \big ( T_{i_m}^{-1}\circ T_1\big) (x)=x.\]
If $n=1$ then $T_1(x)=T_{i_1}(x)$ for  $ x \in Dom(f)$ and some $i_1\geq2$ which leads to a contradiction. Let us assume that $n>1.$
We have
\[ T_{i_1}^{-1} \circ\Big ( O_{m=2}^n \big(T_1 \circ  T_{i_m}^{-1}\big)\Big) \circ T_1(x) =x.\]
  It now yields   that
\[\Big ( O_{m=2}^n \big(T_1 \circ  T_{i_m}^{-1}\big)\Big) \circ T_1(x)=T_{i_1}(x).\]
Let $\theta$ be the function given in assumption (2) of the theorem.
We then have 
 \begin{equation}\label{u0}\theta \circ  \Big ( O_{m=2}^n \big(T_1 \circ  T_{i_m}^{-1}\big)\Big) \circ T_1(x)=\theta\circ T_{i_1}(x)
 \end{equation}  
  We also have that 

\begin{eqnarray*}
\theta \circ  \Big ( O_{m=2}^n \big(T_1 \circ  T_{i_m}^{-1}\big)\Big) \circ T_1(x)
&=&\theta\circ T_1 \circ T^{-1}_{i_2} \circ \Big ( O_{m=3}^n \big(T_1 \circ  T_{i_m}^{-1}\big)\Big) \circ T_1(x)\\
&\geq&\theta\circ T_{i_2} \circ T^{-1}_{i_2} \circ \Big ( O_{m=3}^n \big(T_1 \circ  T_{i_m}^{-1}\big)\Big) \circ T_1(x), \hfill \qquad \qquad \big (\text{Since  } \theta \circ T_1 \geq \theta \circ T_{i_2} \big )\\
&=&\theta \circ \Big ( O_{m=3}^n \big(T_1 \circ  T_{i_m}^{-1}\big)\Big) \circ T_1(x),
\end{eqnarray*}
from which we obtain
\[\theta \circ  \Big ( O_{m=2}^n \big(T_1 \circ  T_{i_m}^{-1}\big)\Big)\circ T_1(x)  \geq \theta \circ \Big ( O_{m=3}^n \big(T_1 \circ  T_{i_m}^{-1}\big)\Big) \circ T_1(x).\]
By repeating the latter argument we obtain that
\begin{equation}\label{u00} \theta \circ  \Big ( O_{m=2}^n \big(T_1 \circ  T_{i_m}^{-1}\big)\Big)\circ T_1(x) \geq \theta\circ T_1(x).\end{equation}
It now follows from (\ref{u0}) and (\ref{u00}) that $\theta\circ T_{i_1}(x) \geq \theta\circ T_1(x)$.
On the other hand by the assumption we have that   $\theta \big(T_1(x)\big)-\theta\big(T_{i_1}(x)\big)$ is non-negative and therefore it must be zero, i.e.,
\[\theta\big(T_1(x)\big)=\theta\big(T_{i_1}(x)\big).\]
By the properties of $\theta$ we  have
$T_1(x)=T_{i_1}(x)$ that contradicts with the fact that $x \in Dom(f).$ This completes the proof of Step (1).\\

{Step 2:} To prove that $f$ and $g$ are measure-aperiodic we need to show that any $T$-invariant probability measure on $\mathcal{B}(X)$ is supported in
$X\setminus D(T)$ where $D(T)=Dom(f)\cap f^{-1} \big ( Dom (g) \big).$  Suppose that $\lambda $ is a probability measure on $\mathcal{B} (X)$ with $T_\# \lambda=\lambda.$
Note first that since $T(x)=x$ for each $x \in X\setminus D(T)$ we have that $ \lambda \big(T^{-1} (A) \big) =\lambda (A)$ for every measurable subset of $D(T).$ It then implies that
$(g \circ f)_\# \lambda =\lambda$ on $D(T).$ Let $f_{|D(T)}$ be the restriction of $f$ on $D(T)$ and let  $\eta$ be the push forward of $\lambda$ by  $f_{|D(T)}.$
Since $(g \circ f)_\# \lambda =\lambda$  on $D(T)$, it  follows that $g_\#\eta=\lambda.$ Let $\mathcal{M}(g, \lambda)$ be the set of positive measures on $\mathcal{B}(X) \cap Dom(g)$
defined by
\[\mathcal{M}(g, \lambda)=\big\{ \zeta ; \, g_\# \zeta=\lambda \big \}.\] Note that $\mathcal{M}(g, \lambda)$ is convex and  $\eta  \in \mathcal{M}(g, \lambda).$
By Theorem \ref{Graf}, extreme points of the set
$\mathcal{M}(g, \lambda)$ are determined by the preimages of $g.$ It follows from the construction  of $g$ that the preimages of $g$ are exactly the  maps  $T_2,..., T_k.$
 Thus, by Theorem \ref{Graf}, the extreme points of  $\mathcal{M}(g, \lambda)$ are exactly the measures, $\zeta_i={T_i}_\# \lambda$ for $i=2,...,k.$ It then follows that
$\eta \in \mathcal{M}(g, \lambda)$ can be written as a convex combination of these measures,
\[\eta= \sum_{i=2}^k\beta_i \zeta_i,\]
where $\beta_2,..., \beta_k$ are non-negative real numbers with $\sum_{i=2}^k\beta_i=1.$ Considering that
 $\eta$ is  the push forward of $\lambda$ by the map $f_{|D(T)}$ together with  $g_\#\eta=\lambda$ we have

 \begin{eqnarray}\label{inter}
 \int_{D(T)} \theta \big (T_1 (x)\big )\, d\lambda=\int_{D(T)} \theta \big (f (x)\big )\, d\lambda \nonumber
&=&\int_{Dom(g)} \theta (y )\, d\eta \nonumber\\
&=&\sum_{i=2}^k\beta_i\int_{Dom(g)} \theta (y )\, d\zeta_i \nonumber\\
&=&\sum_{i=2}^k\beta_i\int_{D(T)} \theta \big (T_i(x)\big )\, d\lambda,
\end{eqnarray}
where $\theta$ is the function given by the assumption (2) in the current Theorem. Since $\theta \big(T_1(x) \big ) \geq \theta \big(T_i(x) \big)$ for $2 \leq i\leq k$, it
follows from (\ref{inter}) that
 \begin{eqnarray*}
 \int_{D(T)} \theta \big (T_1 (x)\big )\, d\lambda
=\sum_{i=2}^k\beta_i\int_{D(T)} \theta \big (T_i(x)\big )\, d\lambda \leq \sum_{i=2}^k\beta_i\int_{D(T)} \theta \big (T_1(x)\big )\, d\lambda=
\int_{D(T)} \theta \big (T_1 (x)\big )\, d\lambda.
\end{eqnarray*}
This in fact implies that \[\sum_{i=2}^k\beta_i\int_{D(T)} \Big [\theta \big (T_i(x)\big )- \theta \big (T_1(x)\big )\Big]\, d\lambda=0.\]
Since each $\beta_i$ is non-negative and $\sum_{i=2}^k\beta_i=1$ at least one of them should be nonzero. Assuming that  $\beta_{i_0}\not=0,$ we must have
 $\theta \big (T_{i_0}(x)\big )=\theta \big (T_1(x)\big )$ for $\lambda$ almost every $x \in D(T).$  Therefore, by the properties of the function $\theta $ we must  have
$T_{i_0}(x)=T_1(x)$ for $\lambda$ almost every $x \in D(T).$ On the other hand, for each $x \in D(T)$ we have $T_{i_0}(x)\not =T_1(x)$ from which we obtain that  $\lambda$ must be zero on $D(T).$
This indeed proves that $\lambda$ must be supported in $X\setminus D(T).$ This completes the proof of Step (2).

  \hfill $\square$\\
 \begin{remark}\label{aper2}
  Theorem \ref{aper} still holds if one replaces the injectivity of $T_2,...,T_k$ with the following assumption,
  \begin{itemize}
   \item
For each $i\geq 2, $   $T_i$ is injective on the set $D_i=\{ x; \, T_1 (x) \not =T_i(x)\}.$
  \end{itemize}
In fact, one just needs to redefine the domain of $f$ and  $g$ as follows:   $Dom(g)= \cup_{i=2}^k T_i(D_i)$ and $Dom(f)=Dom(T_1).$
 \end{remark}

\textbf{Acknowldgement.} I would like to thank Professor Robert McCann for pointing out  a critical issue on the statement of Theorem \ref{main} in the first version of this Manuscript.

\end{document}